\documentclass[a4paper,10pt]{article}
\usepackage[latin1]{inputenc}
\usepackage[T1]{fontenc}
\usepackage{amsmath}
\usepackage{amsfonts}
\usepackage{amstext}
\usepackage{amsthm}
\usepackage[all]{xy}
\usepackage{geometry}
\usepackage{srcltx}
\usepackage{graphics}
\usepackage{epsfig}
\usepackage{subfigure}
\usepackage{slashed}
\usepackage{color}
\usepackage{amssymb}
\usepackage{srcltx}
\newtheorem{theorem}{Theorem}[section]

\newtheorem{defi}[theorem]{Definition}

\title{Not every conjugate point of a semi-Riemannian geodesic is a bifurcation point}
 
\author{Giacomo Marchesi, Alessandro Portaluri$^2$
and Nils Waterstraat}

\begin{document}
\date{}
\maketitle

\footnotetext[1]{{\bf 2010 Mathematics Subject Classification: Primary 58E10; Secondary 53C22, 37J20}}
\footnotetext[2]{The   
author is partially supported by the project ERC Advanced Grant 2013 
No.~339958 ``Complex Patterns for Strongly Interacting Dynamical Systems --- 
COMPAT'', by Prin 2015  ``Variational methods, with applications to 
problems in mathematical physics and geometry'' No.~2015KB9WPT$\_$001 and 
by Ricerca locale 2015  ``Semi-classical trace formulas and their application 
in physical chemistry'' No.~Borr$\_$Rilo$\_$16$\_$01 .} 
\begin{abstract}
\noindent
We revisit an example of a semi-Riemannian geodesic that was discussed by Musso, Pejsachowicz and Portaluri in 2007 to show that not every conjugate point is a bifurcation point. We point out a mistake in their argument, showing that on this geodesic actually every conjugate point is a bifurcation point. Finally, we provide an improved example which yields that the claim in our title is nevertheless true. 
\end{abstract}

\section{Introduction}
Let $(M,g)$ be a semi-Riemannian manifold and $\gamma:[0,1]\rightarrow M$ a geodesic in $M$. A natural question to ask is if there are points on $\gamma$ that can be reached from $\gamma(0)$ by another geodesic, i.e. if there is a geodesic $\widetilde{\gamma}$ starting at $\gamma(0)$ that intersects $\gamma$ at a positive time $t>0$. Of course, here we require that $\gamma$ and $\widetilde{\gamma}$ are \textit{different} in the sense that we exclude the case that $\widetilde{\gamma}(t)=\gamma(s\cdot t)$ for some $s\in[0,1]$. Simple examples are antipodal points on the Riemannian spheres $S^n$ for $n\geq 2$ which are joined by infinitely many distinct geodesics. Note that such antipodal points are conjugate along geodesics, and in general the intuitive meaning of conjugate points suggests that their existence is closely related to the question above.\\
An approach to this problem by variational bifurcation theory was made by the second author with Piccione and Tausk in \cite{PicPorTau}, and with Musso and Pejsachowicz in \cite{MussoPejsachowicz} (cf. also \cite{AleBifIch}). Roughly speaking, a point $\gamma(t^\ast)$ for $t^\ast\in(0,1)$ is called a \textit{bifurcation point} if there is a sequence $\{t_n\}_{n\in\mathbb{N}}$ in $(0,1)$ converging to $t^\ast$ and geodesics $\gamma_n$ starting at $\gamma(0)$, converging to $\gamma$ as $n\rightarrow\infty$ and such that $\gamma_n(t_n)=\gamma(t_n)$ for all $n\in\mathbb{N}$. Bearing in mind that the Jacobi equation along a geodesic is the linearisation of the geodesic equation, it is readily seen from elementary Hilbert manifold theory (cf. \cite{Klingenberg}) and the implicit function theorem that every bifurcation point is a conjugate point. Conversely, it was shown in \cite[Cor. 5.6]{PicPorTau} that every conjugate point in a Riemannian manifold is also a bifurcation point. Moreover, the same is true for causal geodesics in Lorentzian manifolds. The question if this also holds for spacelike geodesics was not considered in \cite{PicPorTau}. Let us recall that there is a crucial difference between causal and spacelike geodesics in Lorentzian manifolds. On causal geodesics there are only finitely many conjugate points as in the Riemannian case. For spacelike geodescics, however, the situation is very different---as Helfer showed in \cite{Helfer}, conjugate points may accumulate. This fact clearly evokes interest in the question of whether along spacelike geodesics likewise every conjugate point is actually a bifurcation point. In \cite{MussoPejsachowicz} an example of a spacelike geodesic in a three dimensional conformally flat semi-Riemannian manifold was discussed that should answer this question in the negative. The aim of this brief note is firstly to point out that there is a mistake in the argument in \cite{MussoPejsachowicz} and that so the question has not yet been answered. Secondly, we use a similar approach to finally construct a geodesic having a conjugate point which is not a bifurcation point.\\
As in \cite{MussoPejsachowicz}, we need to consider the slightly more general setting of perturbed geodesics as in \cite{Marsden}, which we recall in the next section. In particular, we recall a theorem of Jacobi that relates perturbations of geodesics and conformal transformations of metrics. In the third section we revisit at first the two examples from \cite{MussoPejsachowicz}, and for one of them we state two branches of geodesics which show that the considered conjugate point is indeed a bifurcation point. Finally, we use Jacobi's theorem to investigate geodesics on a conformally flat Lorentzian manifold by studying perturbed geodesics in three dimensional Minkowski space.


\section{Geodesics, potentials and bifurcation}
Let $(M,g)$ be a semi-Riemannian manifold of finite dimension $n$ and $V:M\rightarrow\mathbb{R}$ a time-independent potential. Perturbed geodesics (henceforth p-geodesics) connecting two points $p,q\in M$ are critical points of the energy functional

\[E(\gamma)=\frac{1}{2}\int^b_a{g(\gamma'(t),\gamma'(t))\,dt}-\int^b_a{V(\gamma(t))\, dt}\]
that is defined on the Hilbert manifold $\Omega_{pq}$ of all paths $\gamma:[a,b]\rightarrow M$ of regularity $H^{1}$ that connect $p$ and $q$. Alternatively, they are the solutions of the boundary value problem

\begin{equation}\label{geodequ}
\left\{
\begin{aligned}
\frac{D}{dt}\gamma'(t)&+\nabla V(t,\gamma(t))=0\\
\gamma(a)&=p,\quad \gamma(b)=q, 
\end{aligned}
\right.
\end{equation}
where $\nabla V$ denotes the gradient of the potential $V$ with respect to the semi-Riemannian metric $g$. Let us now assume that $\gamma:[a,b]\rightarrow M$ is a p-geodesic of the mechanical system $(M,g,V)$ joining two points $p$ and $q$ in $M$.

\begin{defi}
A point $\gamma(t^\ast)$ for $t^\ast\in(a,b)$ is called a \textit{bifurcation point} if there is a sequence $\{t_n\}_{n\in\mathbb{N}}\subset[a,b]$ converging to $t^\ast$ and p-geodesics $\gamma_n$ of $(M,g,V)$ such that 

\begin{itemize}
\item[(i)] $\gamma_n(a)=\gamma(a)$, $n\in\mathbb{N}$,
\item[(ii)] $\gamma_n(t_n)=\gamma(t_n)$, $n\in\mathbb{N}$,
\item[(iii)] $\gamma'_n(a)\rightarrow \gamma'(a)$ in $T_pM$ as $n\rightarrow\infty$,
\item[(iv)] none of the $\gamma_n$ is a restriction of $\gamma$ to a subinterval of $[a,b]$. 
\end{itemize} 
\end{defi}
\noindent
Let us note firstly that without condition (iv), which was not stated in the definition in  \cite{PicPorTau}, clearly every $\gamma(t)$ would be a bifurcation point. Secondly, that in the unperturbed case condition (iv) is equivalent to requiring that for every $n$, $\gamma_n'(a) \neq s\gamma'(a)$ for some $0<s<1$. Finally, condition (iii) is simply a less technical way to say that $\gamma_n$ converges to $\gamma$ in $\Omega_{pq}$ as required in \cite{MussoPejsachowicz}.\\
The linearisation of the geodesic equation \eqref{geodequ} at the solution $\gamma$ is the \textit{Jacobi equation}

\begin{align}\label{Jacobi}
\frac{D^2}{dt^2}\xi(t)+R(\gamma'(t),\xi(t))\gamma'(t)+D_{\xi(t)}\nabla V(t,\gamma(t))=0,
\end{align}
where $\xi$ is a vector field along $\gamma$, $R$ is the curvature of the connection $D$ and $D_{\xi(t)}\nabla V(t,\gamma(t))$ is the Hessian of $V(t,\cdot)$ with respect to the metric $g$ on the tangent space $T_{\gamma(t)}M$.\\
It is easily seen from elementary Hilbert manifold theory and the implicit function theorem that every bifurcation point $\gamma(t^\ast)$ is a \textit{conjugate point}, i.e. there is a non-trivial solution of the linear boundary value problem

\begin{equation}\label{conjugate}
\left\{
\begin{aligned}
\frac{D^2}{dt^2}\xi(t)&+R(\gamma'(t),\xi(t))\gamma'(t)+D_{\xi(t)}\nabla V(t,\gamma(t))=0\\
\xi(0)&=0,\quad \xi(t^\ast)=0. 
\end{aligned}
\right.
\end{equation}
The following interesting link between geodesics and p-geodesics, that we will need in our construction below, can be found in \cite[Thm. 3.7.7]{Marsden} as \textit{Jacobi's Theorem}.

\begin{theorem}\label{Marsden}
Let $V:M\rightarrow\mathbb{R}$ be bounded above and $c>0$ such that

\[V(p)<c,\quad p\in M.\]
Then the p-geodesics of the mechanical system $(M,g,V)$ with energy $c$ are up to a reparametrisation the same as the geodesics of $(M,g^c)$ with energy $1$ where

\[g^c=(c-V)g.\]
\end{theorem}


\section{Examples -- old and new}
We now discuss the two previous examples from \cite{MussoPejsachowicz}, which are a perturbed and an unperturbed semi-Riemannian geodesic, and explain why they fail to provide examples of conjugate points that are no bifurcation points. Afterwards, we introduce two new examples which indeed show the claim in the title of this paper. 

\subsection{The previous examples}\label{previous}
Let $(M,g_0)$ be the three dimensional semi-Riemannian manifold $\mathbb{R}^3$ with the metric

\[g_0=\frac{d^2}{dx^2}-\frac{d^2}{dy^2}+\frac{d^2}{dz^2}.\] 
Let us consider as in \cite{MussoPejsachowicz} the potential 

\[V(x,y,z)=\frac{1}{2}x^2-\frac{1}{2}y^2+\frac{1}{3}x^3y^3,\]
which has as gradient with respect to $g_0$ the vector field

\[\nabla V(x,y,z)=(x+x^2y^3)\frac{\partial}{\partial x}+(y-x^3y^2)\frac{\partial}{\partial y}\]
and so the equations for a geodesic perturbed by $V$ are

\begin{equation}\label{bvpnonlinI}
\left\{
\begin{aligned}
x''(t)+x(t)+x(t)^2y(t)^3&=0\\
y''(t)+y(t)-x(t)^3y(t)^2&=0\\
z''(t)&=0.
\end{aligned}
\right.
\end{equation}
Obviously, $\gamma_0(t)=(0,0,t)$, $t\in[0,2\pi]$, is a geodesic connecting $(0,0,0)$ and $(0,0,2\pi)$ in $M$, and $(0,0,\pi)$ is a conjugate point along $\gamma_0$. However, the claim in \cite{MussoPejsachowicz}, that there are no non-trivial solutions of \eqref{bvpnonlinI} and that consequently $(0,0,\pi)$ is not a bifurcation point, is wrong.\\
Indeed, the two families of functions 

\[u_\alpha(t)=(\alpha\sin(t),0,t)\quad\text{and}\quad v_\alpha(t)=(0,\alpha\sin(t),t),\quad  t\in[0,2\pi],\, \alpha\in\mathbb{R}\]
are solutions of \eqref{bvpnonlinI}. Hence $(0,0,\pi)$ is a bifurcation point along $\gamma_0$, where we even have two families of solutions branching off from $\gamma_0$.\\
Let us now consider Example 4.7 in \cite{MussoPejsachowicz} which attempts to show that there is a non-perturbed geodesic in $M$ for a metric $g$ such that $(0,0,\pi)$ is a conjugate point but not a bifurcation point of geodesics. Here $g$ is the conformally flat metric $g=e^{2\rho}g_0$, where

\[\rho(x,y,z)=\frac{1}{2}y^2-\frac{1}{2}x^2+\frac{1}{3}x^3y^3.\]   
It is claimed that $\gamma_0(t)=(0,0,\sqrt{2}\,t)$ is a geodesic of energy $1$ having $(0,0,\sqrt{2}\,\pi)$ as a conjugate point, which is a minor inaccuracy. Indeed $\gamma_0$ is a geodesic, however depending on which interval we consider it, it has either an energy greater than $1$ or it does not contain the conjugate point $(0,0,\sqrt{2}\,\pi)$. However, as the neccessary changes to ensure an energy of $1$ would only affect some constants, the following arguments would work verbatim also in that case.\\
Assuming that the energy of $\gamma_0$ is $1$, a geodesic that is close to $\gamma_0$ can be assumed to have energy $1$ as well simply by changing its interval of definition slightly. By Theorem \ref{Marsden}, after a possible reparametrisation, such a geodesic $\gamma(t)=(x(t),y(t),z(t))$ close to $\gamma_0$ needs to to satisfy the p-geodesic equation

\begin{equation}\label{bvpnonlinII}
\left\{
\begin{aligned}
x''(t)+2e^{2\rho(x(t),y(t),z(t))}(x(t)-x(t)^2y(t)^3)&=0\\
y''(t)+2e^{2\rho(x(t),y(t),z(t))}(y(t)+x(t)^3y(t)^2)&=0\\
z''(t)&=0.
\end{aligned}
\right.
\end{equation}
If we multiply the first equation by $y(t)$, the second by $x(t)$, substract the results and integrate from $0$ to $\lambda\in[0,2\pi]$, it follows that every solution of \eqref{bvpnonlinII} such that $x(0)=y(0)=0$, $x(\lambda)=y(\lambda)=0$ satisfies

\begin{align}\label{integral}
\int^\lambda_0{e^{\rho(x(t),y(t),z(t))}(x(t)^2y(t)^4+x(t)^4y(t)^2)\,dt}=0.
\end{align}  
To conclude from this equality as in \cite{MussoPejsachowicz} that $x\equiv y\equiv 0$ is wrong, and so the previous argument is invalid. Indeed, if we set $y\equiv 0$, then \eqref{integral} is satisfied even though $x$ does not need to be identical zero.\\
Let us investigate the case that $y\equiv0$ in \eqref{bvpnonlinII} further. Then the second equation is satisfied trivially, and the first one becomes

\begin{align}\label{singleequ}
x''(t)+2e^{-x(t)^2}x(t)=0.
\end{align}
The second and third author studied in \cite{AleIchDomain}, \cite{AleIchBall} and \cite{AleIchIndef} bifurcation of Dirichlet problems on shrinking domains, and in the special case of ODE's the main theorem reads as follows.

\begin{theorem}\label{bifurcation}
Let $g:[0,b]\times\mathbb{R}\rightarrow\mathbb{R}$ be a smooth function such that $g(t,0)=0$ for all $t\in[0,b]$, and let us consider for $r\in(0,b)$ the family of boundary value problems

\begin{equation}\label{bifnonlin}
\left\{
\begin{aligned}
x''(t)&+g(t,x(t))=0,\quad t\in[0,r]\\
x(0)&=x(r)=0
\end{aligned}.
\right.
\end{equation}
Then there are sequences $\{t_n\}_{n\in\mathbb{N}}$ in $(0,b)$ and $\{x_n\}_{n\in\mathbb{N}}$ in $C^2[0,b]$ such that

\begin{itemize}
	\item[(i)] $t_n\rightarrow t^\ast$ as $n\rightarrow\infty$ for some $t^\ast\in(0,b)$,
	\item[(ii)] $x_n\neq 0$, $n\in\mathbb{N}$, and $x_n$ is a solution of \eqref{bifnonlin} for $r=t_n$,
	\item[(iii)] $x_n\rightarrow 0$ in $C^2[0,b]$ as $n\rightarrow\infty$,
\end{itemize}
if and only if the boundary value problem

\begin{equation}\label{biflin}
\left\{
\begin{aligned}
x''(t)&+f(t)x(t)=0\\
x(0)&=x(t^\ast)=0
\end{aligned}
\right.
\end{equation}
has a non-trivial solution where $f(t)=\frac{\partial g}{\partial x}(t,0)$.
\end{theorem}
\noindent
Now the corresponding boundary value problems \eqref{biflin} for \eqref{singleequ} are

\begin{equation*}
\left\{
\begin{aligned}
x''(t)&+2x(t)=0\\
x(0)&=x(t^\ast)=0
\end{aligned}
\right.
\end{equation*}
which have a non-trivial solution for $t^\ast=\frac{\pi}{\sqrt{2}}$. Hence there are sequences $\{t_n\}_{n\in\mathbb{N}}$ in $(0,2\pi)$ and $\{x_n\}_{n\in\mathbb{N}}$ in $C^2[0,2\pi]$ such that

\begin{itemize}
\item $t_n\rightarrow \frac{\pi}{\sqrt{2}}$ and $x_n\rightarrow 0$ as $n\rightarrow 0$,
\item each $x_n$ is a non-trivial solution of \eqref{singleequ},
\item $x_n(0)=0$ and $x_n(t_n)=0$ for all $n\in\mathbb{N}$.
\end{itemize}
If we now set $\gamma_n(t)=(x_n(t),0,\sqrt{2}\,t)$, then these are p-geodesics in $(M,g_0,-e^{2\rho})$ which start at $\gamma_0(0)$, are different from $\gamma_0$, satisfy $\gamma_n(t_n)=\gamma_0(t_n)$ for $n\in\mathbb{N}$ and converge to $\gamma_0$ as $n\rightarrow\infty$. Hence 

\[\gamma_0(t^\ast)=\gamma_0\left(\frac{\pi}{\sqrt{2}}\right)=(0,0,\pi)\]
is a bifurcation point of p-geodesics on $\gamma_0$. Let us point out, however, that we have not shown that $(0,0,\pi)$ is also a bifurcation point of geodesics for the metric $g$ on $M$. Indeed, in order to apply Theorem \ref{Marsden} again, we would need to know that $E(\gamma_n)=1$ for infinitely many $n\in\mathbb{N}$.



\subsection{The new examples}
The aim is now to construct two examples of geodesics in $M=\mathbb{R}^3$ which have conjugate points that are not bifurcation points. Our examples could be easily adapted to the metric $g_0$ used previously, however, because of the physical motivation we prefer to work with the $3$ dimensional Minkowski metric

\[g_M=-\frac{d^2}{dx^2}+\frac{d^2}{dy^2}+\frac{d^2}{dz^2}.\]
As in the previous section we begin with a perturbed geodesic, where we now consider the potential

\[V(x,y,z)=-\frac{1}{2}x^2+\frac{1}{2}y^2+x^3y+xy^3,\]
which has as gradient with respect to $g_M$ the vector field

\[\nabla V(x,y,z)=(x-y^3-3x^2y)\frac{\partial}{\partial x}+(y+x^3+3y^2x)\frac{\partial}{\partial y}.\]
The equations for a geodesic perturbed by $V$ are

\begin{equation}\label{bvpnonlinIII}
\left\{
\begin{aligned}
x''(t)+x(t)-y(t)^3-3x(t)^2y(t)&=0\\
y''(t)+y(t)+x(t)^3+3y(t)^2x(t)&=0\\
z''(t)&=0,
\end{aligned}
\right.
\end{equation}  
and again $\gamma_0(t)=(0,0,t)$, $t\in[0,2\pi]$, is a geodesic connecting $(0,0,0)$ and $(0,0,2\pi)$ in $M$. Moreover, we see that $(0,0,\pi)$ is a conjugate point along $\gamma_0$, and we now claim that $(0,0,\pi)$ is not a bifurcation point. Indeed, in our case the method from \cite{MussoPejsachowicz} works. We assume that $\gamma(t)=(x(t),y(t),z(t))$, $t\in[0,2\pi]$, is a p-geodesic of $(M,g_M,V)$ such that $\gamma(0)=\gamma_0(0)=(0,0,0)$ and $\gamma(\lambda)=\gamma_0(\lambda)=(0,0,\lambda)$ for some $\lambda\in[0,2\pi]$. If we now multiply the first equation in \eqref{bvpnonlinIII} by $y(t)$, the second by $x(t)$, substract the second from the first and integrate from $0$ to $\lambda$, we obtain

\[\int^\lambda_0{6\,x(t)^2y(t)^2+y(t)^4+x(t)^4\,dt}=0\]
which indeed shows that $x(t)=y(t)=0$ for all $t\in[0,\lambda]$. Hence $\gamma=\gamma_0$ and there is no bifurcation point on $\gamma_0$.\\
Our final aim, and main result of this paper, is the construction of a non-perturbed geodesic in $M$ for a metric $g$ such that $(0,0,\pi)$ is a conjugate point but not a bifurcation point of geodesics. Again we let $g$ be a conformally flat metric $g=e^{2\rho}g_M$, where now the function $\rho:M\rightarrow\mathbb{R}$ is given by

\[\rho(x,y,z)=\frac{1}{2}x^2-\frac{1}{2}y^2+x^3y+xy^3.\]
If we denote by $\nabla^0\rho$ the gradient of $\rho$ with respect to $g_M$, and by $D^0_XY$ the covariant derivative with respect to the Levi-Civita connection with respect to $g_M$, then the Levi-Civita connection of $(M,g)$ is given by

\begin{align}\label{Levi}
D_XY=D^0_XY+d\rho(X)Y+d\rho(Y)X-g_0(X,Y)\nabla^0\rho.
\end{align}
Let us now consider $\gamma_0:[0,2\pi]\rightarrow M$, $\gamma_0(t)=(0,0,\sqrt{\pi}^{-1}t)$, which has energy

\[E(\gamma)=\frac{1}{2}\int^{2\pi}_0{e^{2\rho(0,0,\sqrt{\pi}^{-1}t)}g_M(\gamma'_0(t),\gamma'_0(t))\,dt}=1.\] 
It clearly follows from \eqref{Levi} and the definitions of $\rho$ and $\gamma_0$ that $\gamma_0$ is a geodesic joining $(0,0,0)$ and $(0,0,2\sqrt{\pi})$. Note that the point $(0,0,\pi)$ is on $\gamma_0$, and we now claim that this is a conjugate point.\\
Indeed, bearing in mind that we use as in \cite{MussoPejsachowicz} the sign convention

\[R(X,Y)Z=-D_XD_YZ+D_YD_XZ+D_{[X,Y]}Z\]  
for the curvature, it can be verified by a straightforward computation that on $\gamma_0$

\[R\left(\frac{\partial}{\partial z},\frac{\partial}{\partial x}\right)\frac{\partial}{\partial z}=\frac{\partial}{\partial x}\quad\text{and}\quad R\left(\frac{\partial}{\partial z},\frac{\partial}{\partial y}\right)\frac{\partial}{\partial z}=\frac{\partial}{\partial y}.\]
As

\[\gamma'_0(t)=\frac{1}{\sqrt{\pi}}\,\frac{\partial}{\partial z},\]
a vector field $\xi(t)=x(t)\frac{\partial}{\partial x}+y(t)\frac{\partial}{\partial y}+z(t)\frac{\partial}{\partial z}$ along $\gamma_0$ is a Jacobi field if and only if

\begin{equation*}
\left\{
\begin{aligned}
x''(t)+\frac{1}{\pi}x(t)&=0\\
y''(t)+\frac{1}{\pi}y(t)&=0\\
z''(t)&=0
\end{aligned}
\right.
\end{equation*} 
and so 

\[\gamma_0\left(\pi^\frac{3}{2}\right)=(0,0,\pi)\]
is a conjugate point along $\gamma_0$.\\
We now claim that $(0,0,\pi)$ is not a bifurcation point of geodesics along $\gamma_0$. Indeed, if $\gamma(t)=(x(t),y(t),z(t))$ is a geodesic that is sufficiently close to $\gamma_0$, then we can assume after changing its interval of definition slightly, that $E(\gamma)=1$. Hence by Theorem \ref{Marsden}, $\gamma$ is a geodesic of the mechanical system $(M,g_M,-e^{2\rho})$ and so it satisfies the equations

\begin{equation}\label{bvpnonlinIV}
\left\{
\begin{aligned}
x''(t)+2e^{2\rho(x(t),y(t),z(t))}(-x(t)-y(t)^3-3x(t)^2y(t))&=0\\
y''(t)+2e^{2\rho(x(t),y(t),z(t))}(-y(t)+x(t)^3+3x(t)y^2(t))&=0\\
z''(t)&=0.
\end{aligned}
\right.
\end{equation}
If we multiply the first equation by $y(t)$, the second by $x(t)$, substract the results and integrate from $0$ to $\lambda\in[0,2\pi]$, it follows that every solution of \eqref{bvpnonlinIV} such that $u(0)=v(0)=0$, $u(\lambda)=v(\lambda)=0$ needs to satisfy

\[\int^\lambda_0{e^{2\rho(x(t),y(t),z(t))}(6\,x(t)^2y(t)^2+y(t)^4+x(t)^4)\,dt}=0.\]  
This clearly shows that $x=y=0$ and so $\gamma=\gamma_0$ implying that there is no bifurcation point on $\gamma_0$.\\
The interested reader might ask why the argument that we used at the end of Section \ref{previous} to show that there is a bifurcation point for p-geodesics does not work in the situation of this example. If we set $x$ (resp.\ $y$) equal to zero in \eqref{bvpnonlinIV} then this forces $y$ (resp.\ $x$) to be identically zero and so we cannot use Theorem \ref{bifurcation} as above.

\vspace{1cm}
Giacomo Marchesi\\
School of Mathematics,\\
Statistics \& Actuarial Science\\
University of Kent\\
Sibson Building\\
Parkwood Road\\
Canterbury\\
Kent CT2 7FS\\
UNITED KINGDOM\\
E-mail: gwm7@kent.ac.uk

\vspace{1cm}
Alessandro Portaluri\\
Department of Agriculture, Forest and Food Sciences\\
Universit\`a degli studi di Torino\\
Largo Paolo Braccini, 2\\
10095 Grugliasco (TO)\\
ITALY\\
E-mail: alessandro.portaluri@unito.it\\
Website: aportaluri.wordpress.com

\newpage

Nils Waterstraat\\
School of Mathematics,\\
Statistics \& Actuarial Science\\
University of Kent\\
Sibson Building\\
Parkwood Road\\
Canterbury\\
Kent CT2 7FS\\
UNITED KINGDOM\\
E-mail: n.waterstraat@kent.ac.uk

\end{document}